\documentclass[12pt]{article}
\usepackage{amsmath,amssymb,latexsym}
\usepackage{graphicx}
\usepackage{epstopdf}
\newtheorem{theorem}{Theorem}[section]
\newtheorem{lemma}[theorem]{Lemma}

\newtheorem{remark}[theorem]{Remark}
\newtheorem{corollary}[theorem]{Corollary}
\newtheorem{definition}[theorem]{Definition}
\newtheorem{problem}[theorem]{Problem}

\def\pf{{\bf Proof }}

\begin{document}
\title{On algebraically integrable domains in Euclidean spaces}
%\address{Department of Mathematics, Bar-Ilan University, Ramat-Gan, 52900; Department of Mathematics, Holon Institute of Technology, Holon; Israel }
%\acknowledgement {The author thanks Mikhail Zaidenberg for drawing the author's attention to the subject of this article and stimulating initial discussions.}
%The work of the second author was supported in part by the NSF Grants DMS 9971674 and
%0002195.}
\author{Mark~ Agranovsky}
\maketitle
%\begin{center}
%Bar-Ilan University, Holon Institute of Technology
%\end{center}
%\begin{\dedication}

%\bigskip
%\hspace{5cm}

\begin{center}{{\it To the memory of Yakov Radyno}}
\end{center}
\begin{abstract}

Let $D$ be an infinitely smooth bounded domain $D$ in $\mathbb R^n, n$ is odd.
We prove that if the  volume cut off from the domain by a hyperplane is an algebraic function of the hyperplane, free of real singular points,
then the domain is an ellipsoid. This partially answers a question
of V.I. Arnold: whether odd-dimensional ellipsoids are the only algebraically integrable domains?

\end{abstract}

\footnote{MSC 2010: 44A12, 51M99; keywords:  volume, hyperplane, section, Radon transform, algebraic function, polynomial, ellipsoid. }

\section{Introduction}\label{S:Intro}

\subsection{Formulation of the problem}
Let $D$ be an infinitely smooth  bounded domain in $\mathbb R^n.$   For any affine hyperplane $\{\langle u,x \rangle =t\}, \ u \in \mathbb R^n \setminus 0, t \in \mathbb R,$ denote
$V_D^{\pm}(u,t)$  the volumes of the portions of the domain $D$ on each side from the hyperplane.  Here $\langle,\rangle$ is the inner product in $\mathbb R^n.$
\begin{definition} A bounded domain $D$ in $\mathbb R^n$ is called {\bf algebraically integrable} if the two-valued function $V_D^{\pm}$ is algebraic.
This means that the function $V_D^{\pm}$ can be obtained as a solution of an algebraic equation, i.e., there exists a nonzero  polynomial
$Q(u_1,...,u_n,t,w)$ of $n+2$ variables such that
\begin{equation}\label{E:Q1}
Q(u_1,...,u_n,t,V_D^{\pm}(u_1,...,u_n,t))=0,
\end{equation}
for all $(u,t) \in \mathbb R^n \times \mathbb R$ for which $V_D(u,t)$ is naturally defined, i.e., such that  $D \cap \{\langle u, x \rangle =t\} \neq \emptyset.$
\end{definition}

In the sequel, we will be considering the volume function $V_D$  defined as
\begin{equation}\label{E:cut}
V_D(u,t)=vol_{n}(D \cap \{ \langle u,x \rangle \leq  t\})=\int\limits_{D \cap \{ \langle u,x \rangle \leq  t\}} dx.
\end{equation}

The problem of describing  algebraically integrable domains goes back to Newton \cite{N},\cite{WW}. In connection with Kepler's law in  celestial mechanics, Newton established that there are no algebraically integrable convex bodies (ovals) in $\mathbb R^2$ (for the fact to be true one has to assume the infinite smoothness of the boundaries).
On the other hand, balls and, more generally, ellipsoids in odd-dimensional spaces are algebraically integrable.

V.I. Arnold raised the problem of generalization of Newton's lemma for higher dimensions:
\begin{problem}
{\rm (\cite{Ar},1990-27, 1987-14)} Do there exist algebraically integrable smooth ovaloids different from ellipsoids in $\mathbb R^n$ with odd $n?$
\end{problem}

Ovaloid means in \cite{Ar} a closed hypersurface bounding a convex body.
V. Vassiliev \cite{Vas},\cite{Vas1},\cite{Vas2} proved that there are no algebraically integrable bounded domains (no convexity is required) with $C^{\infty}$ boundaries in $\mathbb R^n$ for even $n.$
In odd dimensions, the question whether ellipsoids are the only algebraically integrable domains, remains unanswered. In this article we partially answer this question in affirmative, under the condition (satisfied for ellipsoids) that the algebraic extension of the function $V_D$ has no real singular points.  No convexity of the domain is a priory assumed.

\subsection{Main notions}

We will be writing the equation of an affine hyperplane in the normalized form $\langle x,\xi \rangle =t, |\xi|=1.$
Correspondingly, the volume function $V_D(\xi,t)$ becomes a function defined on the cylinder $S^{n-1} \times \mathbb R.$

Furthermore, in our considerations, the algebraicity of $V_D$ with respect to $t,$  rather than with respect to the direction variable $\xi,$ will be playing role.

More precisely, the equation $Q(\xi,t,w)=0$ is assumed algebraic in $t$  and  $Q$ belongs to the polynomial ring $C(S^{n-1})[z,w]$ over the coefficient algebra $C(S^{n-1}).$
This ring consists of all polynomials of  two (complex) variables $z,w$ with coefficients - continuous functions on the unit sphere $S^{n-1}.$

Thus, the volume function $V_D(\xi,t)$ satisfies the equation
\begin{equation}\label{E:Qxi}
Q(\xi,t,V_D(\xi,t))=0,
\end{equation}
where
$$Q(\xi,t,w)=\sum_{i=0}^N \sum_{k=0}^{k_i} q_{i,k}(\xi)t^k w^i$$
and the coefficients $ q_{i,k} \in C(S^{n-1}).$

Thus, our condition for the domain $D$ is even weaker than the algebraic integrability since algebraic dependence of  the volume function $V_D(\xi,t)$ is assumed only with respect to $t.$
We will call such domains, i.e. domains $D$ for which (\ref{E:Qxi}) holds, {\bf algebraically $t$-integrable}.

Denote
$$Res_Q(\xi,t)=Res(Q,\frac{\partial Q}{\partial w})(\xi,t)$$
the resultant, with respect to $w,$ of the polynomials $Q(\xi,t,w)$ and its $w$- derivative $\frac{\partial Q(\xi,t,w)}{\partial w}.$
It is equal to
$$Res_Q(\xi,t)=q_N(\xi,t)^{2N-1}\Pi_{i <j}(w_i(\xi,t)-w_j(\xi,t))^2,$$
where $w=w_i(\xi,t)$ are the roots of the algebraic equation $Q(\xi,t,w)=q_0(\xi,t)+...+q_N(\xi,t)w^N=0.$

The resultant is proportional to the discriminant $Disc_Q$ of $Q$ with respect to $w$. Namely,
$$ Res_Q (\xi,t)=q_N(\xi,t) Disc_Q (\xi,t).$$
Since the discriminant is a polynomial of the coefficients of $Q,$  which are continuous functions of  $\xi$ and polynomials of $t,$ the resultant $Res_Q \in C(S^{n-1})[t]:$
$Res_Q(\xi,t)=\sum_{j}^M r_j(\xi)t^j, r_j \in C(S^{n-1}).$
\begin{definition} \label{T:real}
We will call a algebraically $t$-integrable domain $D \subset \mathbb R^n$  {\bf free of real singularities}  if
$Q(\xi,t,V_D(\xi,t))=0,$ where the function $Q(\xi,t,w)$ is continuous in $\xi \in S^{n-1}$ and a polynomial in $t,w$, and
\begin{enumerate}
\item The leading coefficient of $Res_Q(\xi,t)$ satisfies $r_M(\xi) \neq 0, \xi \in S^{n-1}.$
\item $Res_Q(\xi,t) \neq 0$ for all $(\xi,t) \in S^{n-1} \times \mathbb R.$
\end{enumerate}
\end{definition}

\begin{remark}
 The condition 2 says the following. Let $Z$ be the zero variety $Z=Q^{-1}(0) \subset S^{n-1} \times \mathbb C \times \mathbb C$ and $\pi :Z \to \mathbb C, \pi(\xi,z,w)=w-$ the projection to the $w$-plane.  Then the restriction $\pi\vert_{S^{n-1} \times \mathbb R \times \mathbb C}$ is not ramified.
In turn, this means that $V_D(\xi,t)$ has  no real singular points and  extends in $t$ to the whole real line as a real-analytic function.
\end{remark}

Note that if $n$ is even then there is  no infinitely smooth domain $D$ in $\mathbb R^n$ for which the volume function $V_D(\xi,t)$ extends smoothly through the tangent planes to $\partial D.$
Indeed, at the elliptic (having nonzero principal curvatures of the same sign) points $x_0 \in \partial D,$ the volume $V_D(\nu_{x_0},t)$ behaves as $c(t-t_0)^{\frac{n+1}{2}}$ \cite{GG}, where $\nu_{x_0}$ is the normal vector and $t_0$ is the distance of the tangent plane $T_{x_0}(\partial D)$ from the origin.
If $n$ is even then  $\frac{n+1}{2}$ is non-integer and hence the function $V_D(\nu_{x_0},t)$ is not in the class $C^k, k > \frac{n+1}{2},$ at $t=t_0.$

\subsection{Formulation of the main result}
The main result of this article is the following
\begin{theorem}\label{T:main} Let $n \geq 3$ be odd and $D \subset \mathbb R^n$ be a bounded domain with infinitely smooth boundary. Then $D$ is  an algebraically $t$-integrable domain, free of real singularities, if and only if $D$ is an ellipsoid.
\end{theorem}

Let us illustrate  Theorem \ref{T:main} by a partial case of rationally integrable domains.
\begin{corollary} Let $D$ be a bounded domain in $\mathbb R^n, n$ is odd, with infinitely smooth boundary. Suppose that the volume function $V_D(\xi,t)$ is a rational function with respect to $t$, $V_D(\xi,t)=\frac{A(\xi,t)}{B(\xi,t)}, A, B \in C(S^{n-1})[t],$ without real poles: $B(\xi,t) \neq 0,b_k(\xi) \neq 0,\xi \in S^{n-1}, t \in \mathbb R,$ where $b_k(\xi)$ is the leading coefficient of $B(\xi,t).$ Then $D$ is an ellipsoid.
\end{corollary}

 In this case the algebraic equation for $V_D$ is $Q(\xi,t,V_D(\xi,t))=0$, where $Q(\xi,t,w)=B(\xi,t) w-A(\xi,t).$ The polynomial $Q$ has no multiple roots $w$ and the resultant is $Res(\xi,t)=B(\xi,t). $ All the conditions of Theorem \ref{T:main} are fulfilled due to the conditions for the denominator $B(\xi,t).$ The partial case $B(\xi,t)=1$ corresponds to polynomially integrable domains.

\begin{remark} In the original proof of Newton, a contradiction with algebraicity is obtained for family of sections of the planar domain by straight lines rotating around a fixed point (see the article \cite{P}, where an interesting discussion of the Newton's proof in historical aspect is presented ). In our terms, the Newton'a arguments lead to a contradiction with the assumption of algebraic dependence on $\xi.$ In the article \cite{Vas}, the opposite, families of parallel cross-sections are exploited in the proof, i.e., the assumption of the algebraicity with respect to $t$ is essential. Our condition of the $t$-integrability of a domain is just of  the second  kind.
\end{remark}

\subsection{Examples}

\begin{enumerate}
\item
Let $D=B^n$ be the unit ball in $\mathbb R^n.$ Then
$$V_{B^n}(\xi, t)=V_{B^n}(t)=c \int_{-1}^{t}(1-s^2)^{\frac{n-1}{2}}ds.$$ This Abelian integral represents a polynomial $P(\xi,t)$ in $t$ if $n$ is odd, and is a transcendental function if $n$ is even.
The algebraic equation for $w=V_D(\xi,t)$ writes as
$Q(\xi,t,w)=w-P(\xi,t)=0.$
In this case $N=1, q_N(\xi,t)=1, Res_Q(\xi,t)=1$ and therefore $B^n$ is $t$-integrable and free of real singularities.
\item
More generally, consider the case of ellipsoid in the odd-dimensional space, i.e. an image $E=\mathcal A(B_n)$ of the unit ball under a non-degenerate affine transformation. Then
$$V_E(\xi,t)=det \mathcal A \cdot V_{B^n}\Big(\frac{t}{|\mathcal A^{-1}(\xi)|} \Big).$$ Correspondingly, the algebraic equation for the volume function $w=V_E$ is
$Q(\xi,t,w)=w-P \Big(\frac{t}{|\mathcal A^{-1}(\xi)|}\Big).$ Here $Q$ is a polynomial in $t, w$ with the coefficients-continuous functions of $\xi \in S^{n-1}.$ Thus, $E$ is $t$-integrable. The conditions 1,2 in Definition \ref{T:real} are fulfilled, like in the case of the ball.

Thus,  the conditions in Theorem \ref{T:main} for the domain $D$ to be ellipsoid  are necessary.
\end{enumerate}

\section {Preliminaries}\label{S:Prelim}

From now on, we fix a domain $D$ in $\mathbb R^n, n \geq 3$ is odd, satisfying all the conditions of Theorem \ref{T:main}. Without loss of generality, we can assume that $0 \in D.$

The natural domain of definition of the volume function $V_D(u,t),$ defined in (\ref{E:cut}) is
$$ h_{-}(\xi) \leq t \leq h_{+}(\xi),$$
where
$$h_{-}(\xi)=inf_{x \in D}\langle \xi, x \rangle, h_{+}(\xi)=\sup_{x \in D}\langle \xi, x \rangle,$$
are the support functions of the domain $D.$

Since $V_D(-\xi,-t)=vol(D)- V_D(\xi,t),$ we have
\begin{equation}\label{E:odd}
\frac{d^m V_D}{dt^m}(-\xi,-t)=(-1)^{m-1}\frac{d^m V_D}{dt^m}(\xi,t), m \geq 1.
\end{equation}
By the condition $Q(\xi,t,V_D(\xi,t))=0$ is fulfilled in the  domain of definition of $V_D.$ Since $0 \in D,$ this domain includes the set $S^{n-1} \times (-\varepsilon, \varepsilon)$, where $\varepsilon >0$ is sufficiently small.

The algebraic, with respect to $z,w,$ equation $Q(\xi,z,w)=0$ defines $w=w_{\xi}(z)=w(\xi,z)$ as a multi-valued algebraic function  of  $z \in \mathbb C.$
Namely, the equation $Q(\xi,z,w)=0$ has $N=N(\xi,z)$ roots $w=w_j(\xi,z), j=1,...,N.$   The condition for the leading coefficient of the resultant yields that  the degree $N$ does not depend on $\xi$ and $z.$
The function $V_D(\xi,t)$ coincides in its domain of definition with one of the branches $w_j(\xi,t).$

The branch points of the algebraic function $w(\xi,z)$ are the values $(\xi,z)$ at which two or more roots coincide.  The set of branch points is the zero set of the discriminant $Disc_Q(\xi,z)=0.$

Another type of singular points are poles, i.e., the points $(\xi_0,z_0)$ such that  $\lim \limits_{(\xi,z) \to (\xi_0, z_0)} w(\xi,z)=\infty.$
The leading coefficient of the polynomial $Q$ vanishes at the poles: $q_N(\xi_0,z_0)=0$ and hence $Res_Q(\xi,z_0)=0.$

Thus, the condition for $Res_Q(\xi,t)$ of non-having real zeros provides that the algebraic function $w(\xi,z)$ has no singular points, neither poles nor branch points, for $z$ on the real axis, $Im z=0.$
The branches $w_j(\xi, z)$  are locally holomorphic functions of $z$ away from the singular points.

Introduce the {\it parallel section function}
$$A_D(\xi,t)=vol_{n-1}(D \cap \{\langle \xi,x \rangle =t \})=\int\limits_{D \cap \{\langle \xi,x \rangle =t \}}dx.$$

One one hand, $A_D(\xi,t)$ is the $t$-derivative of $V_D(\xi,t):$
$$A_D(\xi,t)=\frac{dV_D}{dt}(\xi,t), h_(\xi) \leq t \leq h_+(\xi).$$
On the other hand, $A_D(\xi,t)$  is the Radon transform of the characteristic function $\chi_D$ of the domain $D:$
$$A_D(\xi,t)=R\chi_D(\xi,t)=\int_{\langle \xi,x \rangle =t}\chi_D(x)dx.$$
The inversion formula for Radon transform yields
$$1=\chi_D(x)=c \int_{|\xi|=1}\frac{d^{n-1}}{dt^{n-1}}A_D(\xi,\langle \xi, x \rangle )d\xi, x \in D,$$
where $d\xi$ is the Lebesgue area measure on $S^{n-1}.$

Let us apply Laplace operator to the both sides of the equality. Then we have $\Delta \chi_D(x)=0, x \in D,$ in the left hand side.  Applying Laplace operator to the right hand side results in twice
differentiating  in $t$ and increasing   by two the order of the $t$-derivative of $A$ under the sign of the integral. Therefore we obtain:
\begin{equation}\label{E:inv}
\int_{S^{n-1}}\frac{d^{n+1}}{dt^{n+1}}A_D(\xi,\langle \xi, x \rangle )d\xi=\int_{S^{n-1}}\frac{d^{n+2}V_D}{dt^{n+2}}(\xi,\langle \xi, x \rangle )d\xi=0, x \in D.
\end{equation}

\section{Proof of Theorem \ref{T:main}}

\subsection{Outline of the proof}
The proof is based on a reduction to the case of polynomially integrable domains, which is solved recently in \cite{KMY},\cite {Ag}.

Assuming that $0 \in D,$ we derive from the inversion formula (\ref{E:inv}) for the Radon transform that the Fourier coefficients  of the decomposition of $V_D(\xi,z),$ where $z$ is near 0,  into spherical harmonics on the unit sphere $\xi \in S^{n-1}$ are polynomials in $z.$

Using the conditions for the resultant we prove that
$V_D(\xi,t)$ has an algebraic continuous extension along two closed paths $\Gamma^{\pm}$ surrounding the singular set in the upper and in lower half-planes, respectively.
Then, using the above polynomiality of the Fourier coefficients  we prove that those continuous branches holomorphically extend, in $z,$ inside the paths.
This removes singularities (branching and poles) and shows that the germ $V_D(\xi,z),$ for $z$ near 0, belongs to a branch which is an entire function of $z \in \mathbb C.$

Great Picard Theorem yields that entire algebraic functions are polynomials. Therefore the germ $V_D(\xi,t)$ extends to the real line as a polynomial in $t$
and we conclude that $D$ is a polynomially integrable domain. According to (\cite{KMY}, \cite{Ag}), $D$ is an ellipsoid.

.

\subsection{Local polynomiality of Fourier coefficients of the volume function $V_D(\xi,t)$ }

Applying translation, we can assume that $0 \in D.$ For sufficiently small $t,$ the hyperplanes $\langle \xi, x \rangle=t$ intersect $D$ for all $\xi
\in S^{n-1}$ and therefore the function $V_D(\xi,t)$ is well defined for $|t|<\varepsilon$ if $\varepsilon$ is sufficiently small.

Decompose $V_D(\xi,t)$ in Fourier series on the sphere $\xi \in S^{n-1}:$
\begin{equation}\label{E:Fourier}
V_D(\xi,t)=\sum_{k=0}^{\infty}\sum_{\alpha=1}^{ d_k} v_{k,\alpha}(t)Y_k^{(\alpha)}(\xi),
\end{equation}
where $\{Y_k^{(\alpha)}\}_{\alpha=1}^{d_k}$ is an orthonormal basic in the space $H_k$ of all spherical harmonics of degree $k.$

The corollary of the inversion formula is the following lemma which plays a key role in the proof of the main result.

\begin{lemma}\label{L:a_k}
The Fourier coefficients $v_{k,\alpha}(t)$ in (\ref{E:Fourier}) are polynomials in the interval $(-\varepsilon,\varepsilon).$
\end{lemma}

\pf

We will be using the notation $V_D^{(m)}(\xi,t)=\frac{d^m V_D}{dt}(\xi,t).$

The algebraic function $w(\xi,z)=V_D(\xi,z),$ and its $z$-derivatives are holomorphic in a small disc $|z|< \varepsilon.$ By Cauchy formula
\begin{equation}\label{E:Cauchy}
V_D^{(n+2)}(\xi,t)=\frac{1}{2\pi i} \int_{|\zeta|=\varepsilon_1}\frac{V_D^{(n+2)}(\xi,\zeta)  d\zeta}{\zeta-t},
\end{equation}
where $0< \varepsilon_1 < \varepsilon, |t|<\varepsilon_1.$

Substituting this representation into (\ref{E:inv}) yields:
$$\frac{1}{2\pi i} \int_{|\xi|=1} \int_{|\zeta|=\varepsilon_1} \frac{V_D^{(n+2)} (\xi,\zeta)} {\zeta- \langle x,\xi \rangle} d\zeta d\xi=0,$$
where $ x \in D$ is in a neighborhood of $0 \in D.$

Let $0< \varepsilon_2 <\varepsilon _1.$ If $|x|<\varepsilon_2$ then $|\langle x,\xi \rangle | < \varepsilon_2$ and the expansion of the Cauchy kernel into the power series  converges uniformly with respect to $\xi.$
Therefore integration the series with respect to $\xi$ is possible
 and one obtains:
$$ \sum_{j=0}^{\infty}\int_{|\xi|=1}b_j(\xi)\langle x,\xi \rangle^j d\xi=0,$$
where
$$ b_j(\xi)= \frac{1}{2\pi i} \int_{|\zeta|=\varepsilon_1} \frac{V_D^{(n+2)}(\xi,\zeta)d\zeta}{\zeta^{j+1}}.$$
Since $n$ is odd, it follows from the relation (\ref{E:odd}) with $m=n+2$ that
\begin{equation}\label{E:bjxi}
b_j(-\xi)=(-1)^j b_j(\xi).
\end{equation}

By homogeneity, each term in the series is zero:
$$ \int_{|\xi|=1} b_j(\xi) \langle x, \xi \rangle ^j d\xi=0, j=0,1,...$$
When $x$ runs over an open $\varepsilon$- neighborhood of $x=0$, the homogeneous polynomials $\psi_{x,j}(\xi)= \langle x, \xi \rangle ^j $ span the space $\mathcal P_j$ of all homogeneous polynomials in $\mathbb R^n$ of degree $j.$
The restriction of this space to the unit sphere $|\xi|=1$ decomposes into the orthogonal sum
$$ span \{ \psi_{x,j}\vert_{S^{n-1}}, |x| < \varepsilon \}=\mathcal P _j\vert_{S^{n-1}} =\oplus_{s=0}^{[\frac{j}{2}]} H_{j-2s},$$
 where $H_k$ is the space of all spherical harmonics of degree $k.$
Thus, we obtain
$$
b_j \perp\oplus_{s=0}^{[\frac{j}{2}]} H_{j-2s} .
$$
On the other hand, the symmetry relation (\ref{E:bjxi}) implies that
$b_j$ decomposes on $S^{n-1}$ into the spherical harmonics of the same parity with $j.$
Therefore, $b_j$ is orthogonal to {\it all} spherical harmonics of degree $k \leq j:$
\begin{equation}\label{E:bj}
b_j \perp\oplus_{k=0}^j H_k.
\end{equation}
Denote $p_{k,\alpha}(t)$ the Fourier coefficients, with respect to $\xi$, of $V_D^{(n+2)}(\xi,t):$
$$ p_{k,\alpha}(t)=\langle V_D^{(n+2)}(\cdot,t), Y_k^{(\alpha)}\rangle_{L^2(S^{n-1})} = \int_{|\xi|=1}V_D^{(n+2)}(\xi,t)Y_{k}^{(\alpha)}(\xi)d\xi.$$
Substitute the integral representation  (\ref{E:Cauchy})  of $V_D^{(n+2)}:$
\begin{equation}\label{E:pk}
p_{k,\alpha}(t)=\frac{1}{2 \pi i} \int_{|\xi|=1}\int_{|\zeta|=\varepsilon_1} \frac{V_D^{(n+2)}(\xi,\zeta)}{\zeta-t}  Y_{k}^{(\alpha)}(\xi)d\zeta d\xi.
\end{equation}
Expansion once again the Cauchy kernel $\frac{1}{\zeta-t}$ into the power series
for $|t| <\varepsilon_2$ yields
$$p_{k,\alpha}(t)=\sum_{j=0}^{\infty}\Big( \int_{|\xi|=1} b_j(\xi)Y_k^{(\alpha)}(\xi)d\xi \Big)t^j.$$
The orthogonality relation (\ref{E:bj}) yields that all the terms with $k \leq j$ vanish and hence $p_{k,\alpha}(t)$ contains only terms with $j<k-1.$ Hence $p_{k,\alpha}$ is a polynomial of degree  $ < k-1.$
Since
$$p_{k,\alpha}(t)=\frac{d^{n+2}v_{k,\alpha}}{dt^{n+2}}(t), $$
the functions $v_{k,\alpha}(t)$ are  polynomials, by successive integration.

\subsection{Removing complex singularities}

Denote $Reg(Q)=\{z \in \mathbb C: Res_Q(\xi,z) \neq 0, \forall \xi \in S^{n-1}\}$ and $Sing(Q)=\mathbb C \setminus Reg(Q).$
\begin{lemma}\label{L:R}
$Sing (Q)$ is a bounded set.
\end{lemma}

\pf

By the condition of Theorem \ref{T:main}, the
leading coefficient $r_M(\xi)$ in the decomposition
$$Res_Q(\xi,t)=\sum_{j=0}^M r_j(\xi)t^j$$ does not vanish for all $\xi \in S^{n-1}.$
Since $r_M \in C(S^{n-1})$ , there exists $c_M >0$ such that  $|r_M(\xi)| \geq c_M, \xi \in S^{n-1}.$ Then for any $z \in \mathbb C$ holds
$$|Res_Q(\xi,z)| \geq c_M|z|^M- c_{M-1}|z|^{M-1}-...-c_0,$$
where $c_j=\|r_j\|_{C(S^{n-1})}.$
Therefore if $|z|>R$ and $R>0$ is sufficiently large then $Res_Q(\xi,z) \neq 0$ for all $\xi \in S^{n-1}$ and hence $Sing (Q) \subset \{|z| \leq R\}.$
Lemma is proved.

Note that the condition of Theorem \ref{T:main} implies that the real line consists of regular points:  $\{Im z=0 \} \subset Reg(Q).$

\begin{lemma} \label{L:Gamma} Let $\Gamma \subset Reg(Q)$ be a Jordan curve , containing a real segment $I_{\varepsilon}= \{ -\varepsilon \leq Re z \leq \varepsilon, Im z=0 \} $ and enclosing the singular set $ Sing(Q).$
Fix a point $z_0 \in \Gamma \setminus I_{\varepsilon}.$
Then the volume function $V_D(\xi,t)$ extends from $(\xi,t) \in S^{n-1} \times [-\varepsilon,\varepsilon]$ to $S^{n-1} \times (\Gamma \setminus \{z_0 \})$ as a continuous function $W_{\Gamma}(\xi,z)$
satisfying $Q(\xi,t,W_{\Gamma} (\xi,z))=0$ for $(\xi,z) \in \ S^{n-1} \times (\Gamma \setminus \{z_0\}).$
\end{lemma}

\pf

Let $z=\varphi(s), s \in [0,1], \varphi(0)=\varphi(1)=z_0, $ be a parametrizations of  the closed curve $\Gamma.$ .
Then the algebraic equation $Q(\xi,z, w)=0$ takes on $\Gamma$ the form
$$Q(\xi,\varphi(s),w)= \sum_{j=0}^{N} q_j(\xi,\varphi(s))w^j=0.$$
The leading coefficient $q_N(\xi,\varphi(s)) \neq 0$ for $(\xi,s) \in S^{n-1} \times [0,1]$ because $Res_Q(\xi,z)= q_N(\xi,z)Disc(Q)(\xi,z) $ and $Res_Q(\xi,z) \neq 0$ when $z \in \Gamma.$

Therefore the equation can be written in the monic form by dividing by the leading coefficient:
\begin{equation}\label{E:f}
f_0(\xi,s)+f_1(\xi, s)w+...+f_{N-1}(\xi,s)w^{N-1}+ w^N=0,
\end{equation}
where
$$f_j(\xi, s)=\frac{q_j(\xi,\varphi(s))}{q_N(\xi,\varphi(s))}.$$
Thus, we deal with an algebraic {\it monic} equation for $w$ with the coefficients-continuous functions on the cylinder $S^{n-1} \times [0,1].$
Since $Res_Q (\xi,z) \neq 0$ on $\Gamma,$ the equation has no multiple roots.

The monodromy theorem \cite{Hu}, Thm. 16.2 (see also \cite{GL}) implies
the algebraic equation (\ref{E:f})  is completely solvable.
This means that there is no monodromy on the cylinder $S^{n-1} \times [0,1]$  and there exist $N$ continuous functions
$W_1(\xi,s),...,W_N(\xi,s)$ on $ S^{n-1} \times [0,1]$ satisfying the equation (\ref{E:f}).

Indeed, consider the mapping $ p:E \to B$
of the space $$E=\{(\lambda_1,...,\lambda_N) \in \mathbb C^N: \lambda_i \neq \lambda_j, i,j=1,...,N.\}$$
to the space $$B=\{(a_0,...,a_{N-1}) \in \mathbb C^{N}: Res(P,\frac{\partial P}{\partial w}) \neq 0 , P(w)=w^N+a_{N-1}w^{N-1}+ \cdots +a_0 \},$$
defined as follows: $p(\lambda)$ is the vector of the coefficients of the monic polynomial of degree $N$ with the roots $\lambda_1,...,\lambda_N.$
By Implicit Function Theorem, $p$ is a $N!$ - covering map. It is regular due to $a_N=1.$
 Since $n \geq 3,$ the fundamental group $\pi(S^{n-1} \times [0,1])=0$ and
by the monodromy theorem the continuous map $f: S^{n-1} \times [0,1], f(\xi,s)=(f_0(\xi,s),\cdots, f_{N-1}(\xi,s))$ can be lifted to a continuous mapping $ W :S^{n-1} \times [0,1] \to E$ such that $pW=f.$ Then $W(\xi,s)=(W_1(\xi,s),\cdots, W_N(\xi,s))$ , where $W_i(\xi,s)$ are the continuous roots of (\ref{E:f}).

Let us regard the functions $W_i(\xi,s)$ as the continuous functions $W_i(\xi,z), z=\varphi(s),$ of
$(\xi, z) \in S^{n-1} \times (\Gamma \setminus \{z_0\}).$  A monodromy can occur at the initial-end point $z_0$ after moving along the closed curves $\Gamma$, so that  one can assign
to $W_i(\xi,z_0)$  two values , corresponding to the values $s=0$ and $s=1$ of the parameter on the curves.

Since the curve $\Gamma$ belongs to the regular set $Reg(Q),$ the functions $W_i(\xi,z)$ are locally holomorphic in a neighborhood $U_{i,\xi,z}$ of  $z \neq z_0.$
Moreover, due to  compactness of $S^{n-1},$ the neighborhood can be chosen the same for all $\xi \in S^{n-1}.$

The germ  $w=V_D(\xi,t), |t|<\varepsilon,$ of the algebraic function $Q(\xi,t, w)=0,$  coincides with one of the functions $W_i, i=1,...,N.$  This function $W_{\Gamma}(\xi,t)$ is just what we need. Lemma is proved.

\begin{lemma}\label{L:ext} The volume function $V_D(\xi,z)$ extends from $|z|<\varepsilon$ to the whole complex plane as an entire function $F(\xi,z)$ of $z \in \mathbb C.$
\end{lemma}

\pf
 Choose in Lemma \ref{L:Gamma} the Jordan curve $\Gamma=\Gamma^+ \subset Reg(Q)$  in the upper half-plane $Im z \geq 0.$   For instance, take
 $$\Gamma^{+}=[-R,R] \cup C_R^{+},$$
 where $C_R^+$ is the open upper half-circle of radius $R$ with the center 0.
The condition of Theorem \ref{T:main} and Lemma \ref{L:R} imply that the curve $\Gamma^{+}$  encloses the singular set $Sing(Q)$ if $R>0$ is sufficiently large. The marked point $z_0$ is taken in the
half-circle $C_R^+$ and the parametrization $\varphi:[0,1] \to \Gamma^+$ satisfies $\varphi(0)=\varphi(1)=z_0.$
Let $W_{\Gamma^+}$ be the branch of $w(\xi,z)$  along $\Gamma \setminus \{z_0\}$  from Lemma \ref{L:Gamma}.

Now fix a basic spherical harmonics $Y_k^{(\alpha)}$ of degree $k$ and consider the functions
$$I_{k,\alpha}(z)=\int_{\xi \in S^{n-1}} Y_k^{\alpha}(\xi) W_{\Gamma^{+}}(\xi, z)d \xi.$$  It is holomorphic in a neighborhood of $\Gamma^+ \setminus \{z_0 \}.$
When $z$ is real, $z=t,$ and  $|z|=|t|<\varepsilon,$ then   $W_{\Gamma^+}(\xi,z)=V_D(\xi,t)$ and therefore
$$I_{k,\alpha}(t)=v_{k,\alpha}(t)$$
is the Fourier coefficient of the volume function $V_D(\xi,t).$

By Lemma \ref{L:a_k}, the Fourier coefficient $v_{k,\alpha}(t)$ is a polynomial in $t$ near $t=0$ (of  degree  depending on $k$).  Therefore, $I_{k,\alpha}(t)$ are polynomials in $t$ in a real
neighborhood of $t=0$ and by analyticity,  $I_{k,\alpha}(z)$ coincides with a complex polynomial $P_{k,\alpha}(z)$ on $\Gamma \setminus \{z_0\}.$ By continuity
it happens at the point $z_0$ as well.

In particular, the one-sided limits
$$I_{k,\alpha}(z_0-0)=\lim\limits_{s \to 1-0}I_{k,\alpha}(\varphi(s)), I_{k,\alpha}(z_0+0)=\lim\limits_{s \to 0+0}I_{k,\alpha}(\varphi(s))$$
 at the point $z_0 \in \Gamma^+$ along the curve $\Gamma^{+}$  coincide:
$$I_{k,\alpha}(z_0-0)=I_{k,\alpha}(z_0+0)=P_{k,\alpha}(z_0).$$

Going back to the definition of $I_{k,\alpha}(z)$ we have:
$$\int_{S^{n-1}} Y_{k}^{(\alpha)}(\xi)W_{\Gamma^+}(\xi,z_0-0)d\xi=\int_{S^{n-1}} Y_{k}^{(\alpha)} (\xi)W_{\Gamma}(\xi,z_0+0)d\xi.$$
Since it is true for all basic spherical harmonics, we conclude that
$$W_{\Gamma^+}(\xi,z_0-0)=W_{\Gamma}(\xi,z_0+0),$$
for all $\xi \in S^{n-1}$ and thus  $W_{\Gamma^+}(\xi,z)$ is continuous on the entire curve $\Gamma^+,$ including the point $z_0.$ .

Furthermore, for  any $m=0,1,...$ one can write
$$\int_{\Gamma^{+}}I_{k}^{(\alpha)}(z) z^m dz=\int_{\Gamma^+} P(z) z^m dz=0.$$
Substituting the expression for $I_{k,\alpha}(z)$ and rewriting the left hand side by Fubini theorem one obtains
$$\int_{S^{n-1}}Y_k^{(\alpha)}(\xi) \Big(\int_{\Gamma^+} W_{\Gamma^+}(\xi,z)z^mdz \Big) d\xi=0.$$

Due to the arbitrariness of the basic spherical harmonic $Y_k^{(\alpha)}$ we conclude
$$\int_{\Gamma}W_{\Gamma^+}(\xi,z)z^mdz=0,$$
for any $\xi \in S^{n-1}.$
Vanishing of all complex moments on $\Gamma$ means that for any $\xi \in S^{n-1}$ the function  $W_{\Gamma^+}(\xi,z)$ is the boundary value of a function
$F_{\xi}^+ (z),$ holomorphic in the domain $\Omega^+$ in the upper half-plane, bounded  by $\Gamma^+.$
Thus, $W_{\Gamma^+}(\xi,z)$ extends from a neighborhood of the curve $\Gamma^+$ inside the domain $\Omega^+$  as a single-valued holomorphic function $F_{\xi}^+(z)$.

Repeating the same argument for the lower half-plane $Im z \leq 0$ and for a corresponding closed curve $\Gamma^{-},$
we conclude that  $V_D(\xi,z)$ continuously extends from the small disc $|z|<\varepsilon$ as a function $F_{\xi}^{-}(z)$ which  is a single-valued holomorphic function of $z$ in the domain $\Omega^{-}$ in the lower half-plane, bounded by $\Gamma^{-}.$

Since $F_{\xi}^+(t)=F_{\xi}^{-}(t)=V_D(\xi,t)$ for $|t|<\varepsilon,$ the two functions define a function $F(\xi,z),$ holomorphic with respect to $z$ in $\Omega=\Omega^+ \cup \Omega^{-}.$

The function $F(\xi,z)$ solves the  equation $Q(\xi,t, F(\xi,t))=0$ for $|t|<\varepsilon$ and, by analyticity, for all $z \in \Omega.$  There is no singularities of the branch $F(\xi,z)$ in the domain $\Omega.$  However, by the construction, the domain $\Omega$ contains {\it all} singular points of the algebraic function $w=w(\xi,z).$
Therefore, the branch $F(\xi,z)$ (as an algebraic function of $z$) has no poles and branch points in $\mathbb C$ and hence is an entire function of $z \in \mathbb C.$
Lemma is proved.

\subsection {End of the proof of Theorem \ref{T:main}}
\begin{lemma} \label{L:Picard} Let $F(z)$ be an entire function in the complex plane. If $F$ is algebraic then $F$ is a polynomial.
\end{lemma}

\pf If $F$ is not a polynomial then $\infty$ is an essential singularity of $G.$ According to Great Picard Theorem, $F(z)$ takes in a neighborhood of $\infty$ any
value, with possibly one exception, infinitely many times.
Therefore, if $F$ is not a polynomial and satisfies a nontrivial algebraic equation $Q(z,F(z))=0, Q$ is a polynomial, then for all $a \in \mathbb C,$ except for at most one, we have
$Q(z_{\nu},a)=Q(z_{\nu}, F(z_{\nu}))=0$ for an infinite sequence of $z_{\nu} \in \mathbb C.$ Since $Q$ is a polynomial ,  $Q(z,a)=0$ for all $z \in \mathbb C$ and for all but one
$a \in \mathbb C.$ Therefore $Q$ is identically zero which is not the case. Lemma is proved.
\begin{corollary} \label{L:polynom} The domain $D$ is polynomially integrable, meaning that for any $\xi \in S^{n-1}$ the volume function $V_D(\xi,t)$ is a polynomial in $t.$
 \end{corollary}
 \pf We have proven in Lemma \ref{L:ext} that the germ $V_D(\xi,z), |z| <\varepsilon,$  extends  to a single-valued branch $F(\xi,z), z \in \mathbb C,$ of the algebraic function $w=w(\xi,z)$ defined by the equation $Q(\xi,z,w)=0.$ By Lemma \ref{L:Picard}, $F(\xi,z)$ is a polynomial of $z.$ Lemma is proved.

To finish the proof, we  refer to the two recent results:
\begin{theorem} \cite{Ag} \label{T:A} Let $D$ be a bounded polynomially integrable domain in $\mathbb R^n, n$ is odd, with infinitely smooth boundary. Then $D$ is convex.
\end{theorem}
\begin{theorem}  \cite{KMY} \label{T:K} Let $D$ be a bounded convex polynomially integrable domain in $\mathbb R^n, n $ is odd, with infinitely smooth boundary. Then $D$ is an ellipsoid.
\end{theorem}

By Corollary \ref {L:polynom}, our domain $D$ is polynomially integrable.  By Theorem \ref{T:A}, $D$ is convex, and by Theorem \ref{T:K} $D$ is an ellipsoid.
Theorem \ref{T:main} is proved.

Bar-Ilan University; Holon Institute of Technology; Israel

{\it E-mail}: agranovs@math.biu.ac.il

\end{document}